\documentclass[11pt,titlepage]{article}
\usepackage[margin=1.2in]{geometry}
\usepackage{amsmath, amsbsy, amsthm, array, mathrsfs}
\usepackage{graphics}
\usepackage{physics}
\usepackage{epsfig, amssymb,latexsym,verbatim}
\usepackage{graphicx}
\usepackage{color}
\usepackage{dsfont, bbm}
\usepackage{relsize}
\usepackage{subcaption}
\usepackage{textcase}

\newcommand{\R}{\mathbb{R}}

\newcommand{\noin}{\noindent}
\newcommand{\bee}{\begin{eqnarray*}}
\newcommand{\ene}{\end{eqnarray*}}
\newcommand{\bec}{\begin{center}}
\newcommand{\enc}{\end{center}}
\newcommand{\be}{\begin{equation}}
\newcommand{\ee}{\end{equation}}

\newcommand{\mb}{\mathbf}
\newcommand{\bs}{\boldsymbol}
\newcommand{\tb}{\textbf}
\newcommand{\pend}{$\blacksquare$}
\newcommand{\vs}{\vskip 3mm}
\newcommand{\bi}{\begin{itemize}}
\newcommand{\ei}{\end{itemize}}

\begin{document}

\begin{titlepage}
\begin{center}
        \vspace*{1cm}
\Large
{  
Non-asymptotic
robustness analysis of regression depth median
 \\[4ex]
}
\vspace{0.4cm}

 {\sc Yijun Zuo}\\[2ex]
         {\small {\em  Department of Statistics and Probability,  Michigan State University} }\\
         {\small East Lansing, MI 48824, USA} \\
         {\small zuo@msu.edu}\\[2ex]

{\small{\today}}
\end{center}
\vspace{0.1cm}
\noin
\large{Abstract}

\vskip 3mm
{\small

\noin
The maximum depth estimator (aka depth median) ($\bs{\beta}^*_{RD}$) induced from regression depth (RD) of Rousseeuw and Hubert (1999) (RH99) is one of the most prevailing estimators in regression.
It possesses outstanding robustness similar to the univariate location counterpart.
Indeed,  $\bs{\beta}^*_{RD}$ can, asymptotically, resist up to $33\%$  contamination without breakdown,
in contrast to the $0\%$ for the traditional (least squares and least absolute deviations) estimators (see Van Aelst and Rousseeuw, 2000) (VAR00)).
The results from VAR00 are pioneering, yet they are limited to regression-symmetric populations (with a strictly positive density) and
 the $\epsilon$-contamination and maximum-bias model.\vs
With a fixed finite-sample size practice, the most prevailing measure of robustness for estimators is the finite-sample breakdown point (FSBP) (Donoho and Huber (1983)).
Despite many attempts made in the literature, only sporadic partial results on FSBP 
 for  $\bs{\beta}^*_{RD}$ were obtained whereas 
 an exact FSBP 
 for  $\bs{\beta}^*_{RD}$ remained open in the last twenty-plus years. 
Furthermore, is the asymptotic breakdown value $1/3$ (the limit of an increasing sequence of finite-sample breakdown values) relevant in the finite-sample practice? (Or what is the difference between the finite-sample and the limit breakdown values?). Such discussions are yet to be given in the literature.
\vs

 This article addresses the above issues,
 revealing an intrinsic connection  between the regression depth of $\bs{\beta}^*_{RD}$ and the newly obtained exact FSBP. It justifies the
 employment of  $\bs{\beta}^*_{RD}$ as a robust alternative to the traditional estimators and demonstrates the necessity and the merit of using the FSBP in finite-sample real practice.
\vs

\bigskip
\par

\noindent{\bf Key words and phrase:} finite-sample breakdown point, regression median, regression depth, maximum depth estimator, robustness.
\bigskip
\par
\noindent {\bf Running title:} Robustness of regression depth median.
}

\end{titlepage}
\section{Introduction}
\vs
The notion of depth in the regression was introduced and investigated two decades ago. The regression depth (RD) of Rousseeuw and Hubert (1999) (RH99) is the most popular example in the literature.
One of the primary advantages of the depth notion in the regression is that it can be utilized  to directly
introduce a median-type maximum depth estimator (aka depth median), which can serve as a robust alternative to the traditional
 least squares and least absolute deviations estimators.
\vs
Robustness of the regression depth induced median ($\bs{\beta}^*_{RD}$) was examined by  Van Aelst and Rousseeuw (2000) (VAR00).
It turns out that the regression median can, asymptotically, resist up to $33\%$  contamination without breakdown,
in contrast to the $0\%$ for the traditional estimators.
The result in VAR00 was established for population (regression-symmetric) distributions
 with a strictly positive density, under the
  maximum-bias framework in the asymptotic sense; the $1/3$ asymptotic result was re-obtained in Van Aelst et al. (2002) as the sample size goes to infinite. It is not directly applicable to the fixed finite-sample size practice.
\vs
In the latter scenario, the most prevailing robustness measure is the finite-sample breakdown point (FSBP), introduced by Donoho and Huber (1983) and popularized and promoted
by Rousseeuw (1984),  Rousseeuw and Leroy (1987) (RL87), and Donoho and Gasko (1992), among others.
The FSBP of  $\bs{\beta}^*_{RD}$ was briefly addressed in RH99. One lower bound was established and the limiting value $1/3$ was listed (later proved in Van Aelst et al. (2002)) for data from regression-symmetric populations with a strictly positive density.
\vs
The lower bound of  FSBP
 in RH99 depends on a generic lower bound of maximum RD value in  RH99. It
is approximately $1/(p+1)$ for a large $n$, and can never approach the asymptotic breakdown value $1/3$ (for $p> 2$).  This implies that the lower bound is not sharp. (A similar lower bound of FSBP was also given in Mizera (2002)(M02)).
\vs
The non-sharpness assertion above is directly verified by an exact FSBP result that was sought in the last twenty-plus years but never obtained until this article.
The exact FSBP  that utilizes the maximum RD value reveals an intrinsic connection between the maximum depth and the FSBP of  $\bs{\beta}^*_{RD}$. The higher the depth value of  $\bs{\beta}^*_{RD}$, the more robust  $\bs{\beta}^*_{RD}$. The new exact FSBP can approach the limiting value of $1/3$.
 Furthermore,  a sharp upper bound of the FSBP for  $\bs{\beta}^*_{RD}$ given pioneeringly in this article indicates that in the finite-sample practice 
  $\bs{\beta}^*_{RD}$ might actually resist much  less than
the asymptotic result of $33\%$  contamination, which renders  $1/3$ inaccurate in the finite-sample practice. Findings here justify (i) the legitimacy of employing  $\bs{\beta}^*_{RD}$ as an alternative to the traditional
estimators and (ii) the necessity and merits of examining
the FSBP for $\bs{\beta}^*_{RD}$.
\vs

Throughout, we are concerned with the FSBP of the 
maximum depth estimator $\bs{\beta}^*_{RD}$ for the 
 regression parameter $\bs{\beta}=(\beta_1, \bs{\beta}^{\top}_2)^{\top}\in \R^p$ ($p\geq 2$) in  the  model:
\begin{eqnarray}
y&=&(1,\mathbf{x}^{\top})\boldsymbol{\beta}+{{e}}, \label{eqn.model}
\end{eqnarray}
where  $\top$ denotes the transpose of a vector, random vector $\mathbf{x}=(x_1,\cdots, x_{p-1})^{\top}$
is in $\R^{p-1}$, 
 and random variables $y$ and ${e}$ are in $\R^1$. $\beta_1$ is the intercept term in the model (\ref{eqn.model}).
\vs

Section 2 briefly reviews the history behind the breakdown point and introduces (i) two versions of the FSBP and (ii) the notion of regression depth and depth-induced median. Section 3 establishes
the exact FSBP (as well as a sharp upper bound) for the regression depth median $\bs{\beta}^*_{RD}$.  A long proof of the major result is deferred to an Appendix.
Section 4 is devoted to the comparison of the sharpness of the lower bounds of the FSBP for  $\bs{\beta}^*_{RD}$ in RH99 with that 
 in this article  and reveals that the latter is sharper than the former.
The article ends with concluding remarks,  including a discussion of the irrelevance of the asymptotic breakdown point in finite samples, supported by substantial empirical evidence.
\vs
\section{Finite sample breakdown point and regression depth}

\subsection{Finite sample breakdown point (FSBP)}

The notion of the breakdown point first appeared in Hodges (1967) and later
was generalized by Hampel (1968, 1971). The finite-sample versions of the
breakdown point, including the {\it addition breakdown point} (ABP) and
{\it replacement breakdown point} (RBP), were introduced by Donoho and
Huber (1983) (DH83). They have become the most prevalent quantitative assessments of
the global robustness of estimators, complementing the assessment of (i) local
robustness of estimating functional captured by the influence function approach
(see Hampel {et al}. (1986) and (ii) global robustness of estimating functional
assessed by the asymptotic breakdown point via the maximum-bias approach (see Hampel {et al}. (1986) and Huber (1981)).
\vs
Stimulating and intriguing discussions on the notion of the breakdown point include Donoho (1982),  Rousseeuw (1984), Rousseeuw and Leroy (1987), Lopuha\"{a} and
Rousseeuw (1991), Maronna and Yohai (1991), Lopuha\"{a} (1992),   Donoho and Gasko
(1992), Tyler (1994), M\"{u}ller (1995), Ghosh and Sengupta (1999), Davies (1987, 1990, 1993),  Davies and Gather (2005), Maronna {et al}. (2006), and Liu {et al}. (2017), among others.
\vs
Some authors favor the ABP in the discussion of the
robustness property of estimators, whereas others prefer the RBP, which
they believe is more simple, realistic and generally more applicable.
Zuo (2001) presented some quantitative relationships between the two versions of the finite-sample breakdown point,
rendering the arguments on the preference (precedence)  between the two versions void in many cases. Nevertheless,
for a given estimator, sometimes one version can be more convenient for the derivation of a desired result. This is especially true for the FSBP of $\bs{\beta}^*_{RD}$ as demonstrated in Proposition 3.2.
\vs
 Throughout,
let $\bs{Z}^{n}=\{\bs{Z}_1, \ldots, \bs{Z}_n\}$ be an uncontaminated sample of
size $n$ in ${\R}^p$, where $\bs{Z}_i=(\bs{x}^{\top}_i, y_i)^{\top}, i\in\{1,\cdots, n\} \subset \mathbb{N}$.

\vs
\noin
\tb{Definition 2.1} [DH83]
The finite-sample \emph{addition breakdown point} (ABP) of
a regression estimator $\bs{T}$ at $\bs{Z}^{n}$ in $\R^p$ is defined as
\be
\mbox{ABP}(\bs{T}, \bs{Z}^{n})=\min_{m \in \mathbb{N}}\{\frac{m}{n+m}:
\sup_{\bs{Y}^{m}}\|\bs{T}(\bs{Z}^{n}+ \bs{Y}^{m})-\bs{T}(\bs{Z}^{n})\|=\infty\},
\ee
where $\bs{Y}^m$ denotes a dataset of size $m$ with arbitrary values in $\R^{p}$ and
$\bs{Z}^{n}+ \bs{Y}^m$ denotes the contaminated sample by adjoining $\bs{Y}^m$ to
$\bs{Z}^n$ (i.e., $\{\bs{Z}_1, \cdots, \bs{Z}_n,\bs{Y}_1, \cdots, \bs{Y}_m\}$), $\|\cdot\|$ stands for Euclidean norm.
\vs
\noin
\tb{Definition 2.2} [DH83]
The finite-sample \emph{replacement breakdown point} (RBP) of a regression estimator $\bs{T}$ at $\bs{Z}^{n}$ in $\R^p$
 is defined  as
\begin{equation}
\text{RBP}(\bs{T},\bs{Z}^{n}) = \min_{m \in \mathbb{N}}\bigg\{\frac{m}{n}: \sup_{\bs{Z}_m^{n}}\|\bs{T}(\bs{Z}_m^{n})- \bs{T}(\bs{Z}^{n})\| =\infty\bigg\},
\end{equation}
where $\bs{Z}_m^{n}$
denotes an arbitrary contaminated sample by replacing $m$ original sample points in $\bs{Z}^{n}$ with  arbitrary $m$ points in $\R^{p}$.
\vs
In other words, the ABP and RBP of an estimator are respectively the
minimum addition fraction and replacement fraction of the contamination which could drive
the estimator beyond any bound.
\vs
\subsection{Regression depth of Rousseeuw and Hubert (1999)}

The regression depth (RD) in RH99 was defined based on a notion of \emph{nonfit}. Equivalent definitions  
 were given in the literature, see e.g., VAR00, M02, and Zuo (2020, 2021) (Z20, Z21), also see (\ref{rd-proof.eqn}) below.
A formal definition in Rousseeuw and Struyf (2004) (RS04) is
\be
\mbox{RD}(\bs{\beta}; \mbox{Pr})= \inf_{D\in{\mathcal{D}}}\left\{\mbox{Pr}\left((r(\bs{\beta})\geq 0) \cap D \right)+\mbox{Pr}\left((r(\bs{\beta})\leq 0) \cap D^c \right)\right\},\label{RS04.eqn}
\ee
where $\mbox{Pr}$ is the joint probability distribution of $(\mb{x}^{\top}, y)$ in (\ref{eqn.model}), $\mathcal{D}$ is the set of all vertical closed halfspaces D whose boundary $\partial D$ is parallel to the $y$-axis and $\mbox{Pr}(\partial D)=0$, the complement of $D$ is $D^c$ 
and $r(\bs{\beta}):=y-(1,\mb{x}^{\top})\bs{\beta}$.
The maximum regression depth functional $\mb{T}^*_{RD}$ 
 is defined as (see RH99)
\be \mb{T}^*_{RD}(\mbox{Pr})=\arg\!\max_{\bs{\beta}\in\R^p}\mbox{RD}(\bs{\beta};\mbox{Pr}). \label{T-RD.eqn}
\ee
\vs
 One obtains the sample version of RD$(\bs{\beta};\mbox{Pr})$ and $\mb{T}^*_{RD}(\mbox{Pr})$ by replacing $\mbox{Pr}$ with $\mbox{Pr}_n$, the latter is the empirical distribution based on  a given sample $\mb{Z}^n=\{(\mb{x}^{\top}_i, y_i)^{\top}, i\in \{1,\cdots, n\}\}$ in $\R^{p}$.
(In the empirical case, the RD discussed originally in RH99 divided by $n$ is identical to (\ref{RS04.eqn}) and definition 2.3 below). Hereafter $\mb{Z}^{n}$ and $\mbox{Pr}_n$ will be used interchangeably.
\vs
For examples and illustrations and explanations of the $\mbox{RD}(\bs{\beta};\mbox{Pr}_n)$ defined above in $\R^2$, we refer to RH99 (also see Section 3).
 If there are several $\bs{\beta}$s that attain the maximum depth value on the right-hand side (RHS) of (\ref{T-RD.eqn}), then the average of all those $\bs{\beta}$s is taken.\vs
The following equivalent definition is useful in our proof in the sequel (for explanations of related terms, see Section 3).
\vs
\noin
\tb{Definition 2.3}~
{For any $\bs{\beta}\in \R^p$, 
RD$(\bs{\beta};\mbox{Pr})$
is the minimum probability mass that needs to be passed when tilting (the hyperplane induced from) $\bs{\beta}$ in any way until it is vertical.}\vs

The notions of breakdown point and regression depth seem unrelated and have nothing to do with each other.  But in the next section,  it is shown that they are actually closely connected in the case of $\mb{T}^*_{RD}(\mbox{Pr}_n)$.\vs Note that $\bs{\beta}^*_{RD}(\mbox{Pr}_n)$ also denotes the empirical maximum depth estimator (or depth median). But the latter is defined to have maximum depth whereas this is not necessarily the case for $\mb{T}^*_{RD}(\mbox{Pr}_n)$ (see the proof of
proposition 3.2) which leads us to keep both notations.

\vs
\section{Finite sample breakdown point of regression depth median}

\subsection{A preliminary lemma}\vs
To facilitate discussions and proofs in the sequel, we first provide a characterization of RD in Definition 2.3.
\vs
For a given $\bs{\beta}=(\beta_1,\bs{\beta}^{\top}_2)^{\top}\in\R^p$, denoted (hereafter) by $H_{\bs{\beta}}$  the unique hyperplane determined by $y= (1,\mb{x}^{\top})\bs{\beta}$.
Denote the angle between the hyperplane $H_{\bs{\beta}}$ 
 and the horizontal hyperplane plane $H_h$ (determined by $y=0$) by ${\theta}_{\bs{\beta}}$ (hereafter consider only the acute one).
 That is, ${\theta}_{\bs{\beta}}$ is the angle between the normal vector $(-\bs{\beta}^{\top}_2, 1)^{\top}$ of $H_{\bs{\beta}}$ and the normal vector $(\mb{0}^{\top}, 1)^{\top}$ of $H_h$ in the $(\mb{x}^{\top}, y)^{\top}$-space.  Hence, $$\cos({\theta}_{\bs{\beta}})=\frac{1}{\sqrt{\|\bs{\beta}_2\|^2+1}}.$$
 Therefore, it is not hard to see that $|\tan(\theta_{\bs{\beta}})|=\|\bs{\beta}_2\|$.
\vs

\emph{Tilting} ${\bs{\beta}}$ to a vertical position (some vertical hyperplane $H_v$)  in Definition 2.3 means tilting $H_{\bs{\beta}}$ along a \emph{hyperline} ${l}_v(\bs{\beta})$ (the intersection/common part of  $H_{\bs{\beta}}$ and  $H_v$) to $H_v$ (for illustration in $\R^2$, see Figure 1 of RH99, where a hyperplane reduces to a line and a hyperline becomes an intersection point). The probability mass that needs to be passed
(when tilting $\bs{\beta}$ in any way until it is vertical) is either ${\mbox{Pr}_n\left((r(\bs{\beta})\geq 0)\cap D(H_v))+\mbox{Pr}_n((r(\bs{\beta})\leq 0)\cap D^c(H_v)\right)}$ or $\mbox{Pr}_n((r(\bs{\beta})\geq 0)\cap D^c(H_v))+\mbox{Pr}_n((r(\bs{\beta})\leq 0)\cap D(H_v))$ (see (\ref{RS04.eqn})),  where $D(H_v)$ is one of the closed halfspaces with $H_v$ as its boundary.
The two quantities  represent the two fractions of data points passed by tilting $H_{\bs{\beta}}$ in two ways (clockwise or counter-clockwise, see Figure 1), respectively.
Denote the minimum of them  by
$ \min_{fr}(l_v(\bs{\beta}), \mbox{Pr}_n)$. 
\vs
 Namely, $ \min_{fr}(l_v(\bs{\beta}), \mbox{Pr}_n)$ is  the minimum of the two fractions of data points touched by tilting $H_{\bs{\beta}}$ in the definition of $\mbox{RD}$
 to a vertical position $H_v$ along $l_v(\bs{\beta})$ in two ways (one way is by crossing the double wedge formed by two single wedges with an acute angle between $H_{\bs{\beta}}$ and $H_v$ (the two shaded regions in Figure 1) and the other way is by
 passing through the double wedge formed by two single wedges with an obtuse angle between $H_{\bs{\beta}}$ and $H_v$ (the other two regions in Figure 1).
\vs

\bec
\begin{figure}[h!]
\vspace*{-8mm}
\includegraphics[width=\textwidth]{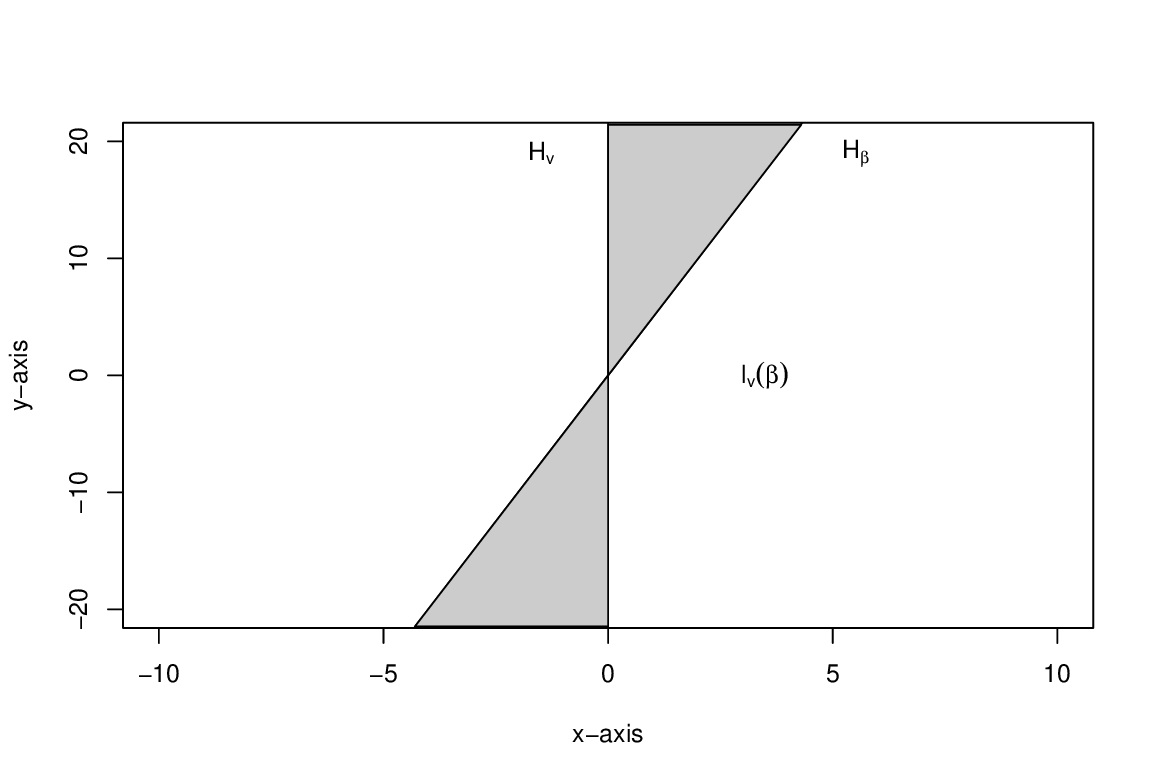}
\vspace*{-13mm}
 \caption{\small A two-dimensional vertical cross-section of a figure in  $\R^p$.  There are two ways to tilt $H_{\bs{\beta}}$ to a vertical position $H_v$ (which does not necessarily contain the origin in the definition, it does in this figure) along hyperline $l_v(\bs{\beta})$ (which passes
 $(0, 2)$ in the figure). One way is crossing the two wedges each with an acute angle (the shaded double wedge), the other way is passing through the other two wedges each with an obtuse angle (the unshaded double wedge). That is, counter-clockwise or clockwise, tilting $H_{\bs{\beta}}$ to $H_v$ along hyperline $l_v(\bs{\beta})$ (a point in the two-dimensional cress-section).
  }
 \label{fig-4-proof}
\end{figure}
\enc

\vspace*{-10mm}
 In other words, within the two-dimensional plane that is perpendicular to the horizontal hyperplane $H_h$ (a vertical cross-section),  tilting $H_{\bs{\beta}}$ can be seen in this two-dimensional plane  in  either a clockwise or counter-clockwise manner.
 \emph{Hereafter, when we talk about clockwise or counter-clockwise tilting, all are in this sense (within the vertical cross-section)}. Above discussions lead to the following equivalent, 
  which is even more useful in the sequel.
 \vs
\noindent
\tb{Lemma 3.1}  For a given data set $\mb{Z}^n$, the regression depth of $\bs{\beta}$ defined in Definition 2.3 can be characterized as
  \be
 \mbox{RD}(\bs{\beta};\mbox{Pr}_n)=\inf_{l_v({\bs{\beta}})}\min_{fr}(l_v(\bs{\beta}), \mbox{Pr}_n), \label{rd-proof.eqn}
 \ee
where the infimum is taken over all possible $l_v({\bs{\beta}})$s or equivalently all possible $H_v$s.
\vs
\noin
\tb{Proof:} ~~
In light of the discussions before the lemma, the proof is trivial. \hfill \pend
\vs
 For a given sample $\mb{Z}^n$ in $\R^p$, 
  write
 \begin{align}
k^*(\mbox{Pr}_n)&=\max_{\bs{\beta}\in\R^p} n\mbox{RD}(\bs{\beta}, \mbox{Pr}_n). \label{k*.eqn}
\end{align}


\noindent
\tb{Remark 3.1}
\bi
\item[]The quantity $n\mbox{RD}(\bs{\beta}, \mbox{Pr}_n)$ in (\ref{k*.eqn}) is the least number of data points touched by tilting $H_{\bs{\beta}}$ to a vertical position  in light of
(\ref{rd-proof.eqn}).
The $k^*(\mbox{Pr}_n)$ then is the maximum (w.r.t. $\bs{\beta}$s) of the least number of data points touched by tilting a $H_{\bs{\beta}}$  to a vertical position. That is the least number of data points touched by tilting $H_{\bs{\beta}_m}$ in any way to a vertical position with $\bs{\beta}_m$ attaining the maximum RD (i.e. $\bs{\beta}_m$ is an RD maximizer).
\hfill \pend
\ei
\vs
\subsection{Upper bounds of finite sample breakdown point for  $\bs{T}^*_{RD}$}

A regression estimator $\mb{T}$ is called \emph{regression equivariant} (page 116 of RL87) if \be \mb{T}\left(\{(\mb{x}^{\top}_i, y_i+(1,\mb{x}^{\top}_i) \mb{b})^{\top}\}\right)=
\mb{T}\left(\{(\mb{x}^{\top}_i, y_i)^{\top}\}\right)+\mb{b}, ~\forall ~i\in\{1,\cdots, n\},~ \mb{b}\in\R^p
\label{regression.equi}
\ee
\vs
It is seen that $T^*_{RD}(\mbox{Pr}_n)$ is regression equivariant (see Z21).  In the sequel, for simlicity we write $T^*_{RD}$ for $T^*_{RD}(\mbox{Pr}_n)$. 
\vs
Based on Theorem 4 of RL 87, one obtains
an RBP upper bound for $\mb{T}^*_{RD}$ as $$~~ \mbox{RBP}(\mb{T}^*_{RD}, \mb{Z}^n)\leq(\lfloor (n-p)/2\rfloor+1)/n,$$ for a given $\mb{Z}^n$ in $\R^{p}$, where $\lfloor \cdot\rfloor$ is the floor function.
A similar ABP bound for $\mb{T}^*_{RD}$ is given below.
\vs
\noindent
\tb{Proposition 3.1}. For a given $\mb{Z}^n$ in $\R^{p}$, we have
\be
~~ \mbox{ABP}(\mb{T}^*_{RD}, \mb{Z}^n)\leq \frac{n-p+1}{2n-p+1}.\label{upper.eqn}
\ee
\vs
\noin
\tb{Proof:} Utilizing part (a) of the Theorem 2.1 of Zuo (2001), one directly obtains the upper bound based on its RBP upper bound $(\lfloor (n-p)/2\rfloor+1)/n$ with $\beta=-p$ and $m=\epsilon=1$.
(Alternatively, one can also prove the result analogously as the proof of Theorem 4 of RL87).\hfill \pend
\vs
The ABP upper bound  for $\mb{T}^*_{RD}$ above, albeit never appearing in the literature, is not sharp. In the following we will establish  sharper upper bounds for both ABP and RBP.\vs
\subsection{The exact ABP and a sharper RBP upper bound for  $\bs{T}^*_{RD}$}

We shall say  $\mb{Z}^{n}$  is \emph{in general position} (IGP)
when any $p$ of observations in $\mb{Z}^{n}$ give a unique determination of $\bs{\beta}$ (i.e., $y_{i_j}=(1,\mb{x}^{\top}_{i_j})\bs{\beta}, i_j\in\{1,2, \cdots, n\}, j\in\{1,\cdots, p\}$ yields a unique solution of $\bs{\beta}$).
In other words, any $(p-1)$ dimensional affine subspace of the space $(\mb{x}^{\top}, y)$ contains at most p observations of
$\mb{Z}^{n}$.
When the observations come from continuous distributions, the event ($\mb{Z}^{n}$ being in general position) happens with probability one.
\vs

\noin
\tb{Proposition 3.2} For a given IGP $\mb{Z}^n$ (or $\mbox{Pr}_n$)  in $\R^{p}$, we have
\begin{align}
(\mbox{A})~~ \mbox{ABP}(\mb{T}^*_{RD}, \mb{Z}^n)&= \frac{k^*(\mb{Z}^n)-p+1}{n+k^*(\mb{Z}^n)-p+1}, \label{lower.eqn}\\[2ex]
(\mbox{B})~~ \mbox{RBP}(\mb{T}^*_{RD}, \mb{Z}^n)&\leq \frac{k^*(\mb{Z}^n)-p+1}{n}.\label{upperRBP.eqn}
\end{align}
\vs
\noindent
\tb{Proof:} see the Appendix. \hfill \pend
\vs
\noindent
\tb{Remarks 3.2}~
\bi
\item[(1)] The assumption that $\mb{Z}^n$ is IGP is crucial in our proof.
It implies that 
any $p$ points from $\mb{Z}^n$ uniquely determine a hyperplane and
any $(p-2)$ ``hyperline" in the $\mb{x}$-space contains at most $p-1$ points with $\bs{x}_i \in \{\bs{x}_i, i\in\{1,\cdots, n\}\}$. If we introduce a quantity $c(\mb{Z}^n)$ which is the maximum number of data points contained in a $p-1$ dimensional hyperplane in $(\mb{x}^{\top}, y)$-space (see Zuo (2019), page 1186, and M02, page 1695), then we can replace the $p$s in Proposition 3.2 with this $c(\mb{Z}^n)$, and the proof can be modified accordingly.
\item[(2)]
 After the establishment of the lower bound in (A), this author learned that
    the same lower bound appeared in Van Aelst et al. (2002) in their proof of the limiting breakdown value $1/3$. There are some essential differences between the two: (i)an \emph{extra} assumption that samples are from regression-symmetric distribution (see RS04 for definition) with a positive density is required there, (ii) their result holds only almost surely, and (iii) our proof is methodologically different from theirs.
 \item[(3)]
  The upper bound of ABP in (A) was never established in the last twenty-plus years whereas the RBP upper bound in (B) is sharper than any existing ones before.
\item[(4)]
Lower bounds of  RBP for  $T^*_{RD}$ exist in the literature, including the ones in Van Aelst et al. (2002), M02, and RH99. The first one holds almost surely and with an extra assumption while the latter two are not sharp.
RH99 presented in their Corollary of Conjecture 1
 an RBP lower bound  for  $T^*_{RD}$  under some assumptions as follows, 
\be\mbox{RBP}(T^*_{RD}, \mb{Z}^n)\geq\frac{1}{n}\left(\Big\lceil \frac{n}{p+1}\Big\rceil-p+1\right),\label{lower bound-rh999}\ee
 which is not helpful to compare it directly with the ABP lower bound  
 in (A). That is, it is not clear which one is sharper at this moment. We tackle this issue next. \hfill \pend
\ei

\noindent
\section{Sharpness of lower bounds of FSBP for $\mb{T}^*_{RD}$}
One naturally wonders given that a lower bound is already established in RH99, what is the merit to have Proposition 3.2?
Or rather is there any difference between the two breakdown values?
 Let us focus on the lower bound of the RBP in RH99. The essential difference between this and the lower bound of the ABP (which is the same as the upper bound of the ABP) in Proposition 3.2 is: the former (RH99) employing a lower bound, $\lfloor n/(p+1)\rfloor$, for the maximum depth value whereas the latter (Proposition 3.2) utilizing the maximum depth value, $k^*(\mb{Z}^n)$, directly. Note that $\lfloor n/(p+1)\rfloor \le k^*(\mb{Z}^n)$ (conjecture 1 of RH99). \vs
 Consequently,
(i) the former, the RHS of (\ref{lower bound-rh999}), purely depends on $n$ and $p$ (this could be an advantage) whereas the latter, $(k^*(\mb{Z}^n)-p+1)/(n+k^*(\mb{Z}^n)-p+1)$,  depends on the configuration of $\mb{Z}^n$. (Note that FSBPs dependent on the configuration of data are not rare, e.g., that of median absolute deviations (MAD), also see Huber (1984) Theorem 3.1 and
Davies and Gather (2007));
(ii) for a fixed $p>2$, as $n\to \infty$ the former approaches $1/(p+1)$ which can never be $1/3$ (the asymptotic breakdown value) while the latter
can approach  $1/3$ with data from the regression-symmetric  population (see RS04 for definition);
(iii) when $p>2$, 
$\mb{T}^*_{RD}$ can resist a much higher fraction of contamination than the one given in the RHS of (\ref{lower bound-rh999}) (see Table 1 below), implying the lower bound in RH99 is not sharp.
\vs

ABP and RBP, albeit employing different contamination schemes, are the fraction (or percentage) of contamination that can force an estimator beyond any bound (becoming useless). From the contamination fraction/percentage interpretation, ABP is usually slightly smaller than RBP for the same estimator. For example, in the sample mean case, it is $1/(n+1)$ v.s. $1/n$; in the sample median case, it is $n/(n+n)$ v.s. $\lfloor(n+1)/2\rfloor/n$.
Directly comparing  ABP and RBP, in terms of their magnitude,  is unfavorable (or unfair) to ABP. But, without a better measure, it is at least one approach to check the sharpness of the two lower bounds.
\vs
To better appreciate the sharpness of the lower bound of the ABP in  Proposition 3.2 and know better the quantitative difference between the two lower bounds of FSBP
 for the same $\mb{T}^*_{RD}$, 
 we carry out a small scale simulation study to
calculate the average differences of  $(k^*(\mb{Z}^n)-p+1)/(n+k^*(\mb{Z}^n)-p+1)$ (the lower bound of the ABP in Prop. 3.2) with 
$\big(\lceil {n}/{(p+1)} \rceil-p+1\big)/n$ (the lower bound of the RBP in RH99) in $1000$ multivariate $N(\mb{0}, \mb{I})$ 
 samples for different small $ns$ and $ps$; the results  are given in  Table \ref{t1}.
\vspace*{-3mm}
\bec
\begin{table}[h!]
\centering
\caption{Average differences (here and after {expressed in percentage points})  between the lower bound of the ABP in Prop. 3.2 and that of the RBP in RH99, based on 1000 standard normal samples.}\vs
\begin{tabular}{c c c c c c c} 
n& 10& 20& 30& 50& 100&200\\[-0.0ex]
\hline\\[.1ex]
p=2& -3.725 &-1.776 &-0.913 &-2.237 &-2.456 &-1.805 \\[.5ex]
p=3& 10.38 & 8.736 &5.235 & 4.646 &5.139 & 5.328\\[.5ex]
p=5& 31.52 & 16.85 &15.09 &12.02 & 11.47 & 11.15 \\[.5ex]
\hline
\end{tabular}
\label{t1}
\end{table}
\vspace*{-4mm}
\enc

\indent
It is readily apparent from the table that the lower bound of the ABP in  Proposition 3.2 is
sharper than the lower bound of the RBP in RH99 because of all the positive entries (when $p>2$) and as well as the negative ones (when $p=2$) since all entries should be negative if the number of contaminating points m is the same in the two contamination schemes.\vs

The positive entries in the table 
imply that  $\mb{T}^*_{RD}$ can resist much higher contamination percentages than what is provided by the lower bound of the RBP in RH99.
For $p=5$, the difference in the table 
 decreases when $n$ gets larger, this is not the case for $p=2$ (i.e., it is not monotonic, e.g., when $p=2$ and $n=500$, the entry in the table will be $-1.352$). For a fixed $n$, when $p$ increases so does the difference.
The lower bound of the RBP in RH99 becomes negative and uninformative if $\lceil {n}/{(p+1)} \rceil <p-1$.  (this explains the unusually large entry in the case $p=5, n=10$).
\vs
The results in the table demonstrate the merit of the lower bound of the ABP in  Proposition 3.2.
One question that might be raised for the results in the table is: are those results distribution-free?
That is, if the underlying distribution of the samples changes, does the lower bound of the ABP in Proposition 3.2 still have any advantage over the one in RH99? To answer this question, we carried out a small scale simulation study and results are reported in Table \ref{t22}.
\bec
\begin{table}[h!]
\centering
\caption{Average differences  between the lower bound of the ABP in Prop. 3.2 and that of the RBP in RH99, based on 1000 contaminated normal samples.}\vs
\begin{tabular}{c c c c c c c} 
n& 10& 20& 30& 50& 100&200\\[-0.0ex]
\hline\\[.1ex]
p=2& -3.687 &-1.089 &-0.929 &-2.261 &-2.604&-2.042 \\[.5ex]
p=3& 10.45 & 8.742 &5.214 & 4.572 & 4.888 & 4.993\\[.5ex]
p=5& 31.36& 16.92 &16.09 &11.85 & 11.30 & 10.78 \\[.5ex]
\hline
\end{tabular}
\label{t22}
\end{table}
\vspace*{-10mm}
\enc
\vs

Here we generated $1000$ samples $\mb{Z}^{(n)}=\{(\mb{x}^{\top}_i, y_i)^{\top}, i\in\{1,\cdots, n\}, \mb{x}_i \in \R^{p-1}\}$ from the Gaussian distribution with the zero mean vector and $1$ to $p$ as its diagonal entries of the diagonal covariance matrix  for various $n$s and $p$s. Each sample is contaminated by $5\%$ i.i.d. normal $p$-dimensional points with individual mean $10$ and variance $0.1$. Thus, we no longer have symmetric errors and a homoscedastic variance model.
\vs
Comparing the table entries in Tables 1 and 2, we conclude that the sharpness of the lower bound of the ABP in Proposition 3.2 over the lower bound of the RBP in RH99  almost does not depend on the underlying distributions overall (this is confirmed in the multivariate t-distribution case). However, $k^*(\mb{Z}^n)$ does depend on the configuration of points of $\mb{Z}^n$.
\vs
\noindent
\section{Concluding remarks}
\vs
\noin
{(I)} \tb{The state of the art on the finite-sample breakdown point of $\mb{T}^*_{RD}$}.\vs
\noin
In the last twenty-three years, despite numerous attempts made, only sporadic partial results on FSBP 
 of $\mb{T}^*_{RD}$ were obtained.
No existing results in the literature offer an exact FSBP or sharp upper bounds as Proposition 3.2 does.
There are several RBP lower bounds in the literature, e.g., in Van Aelst et al. (2002), M02, and RH99 for  $\mb{T}^*_{RD}$. But the one in Van Aelst et al. (2002) holds almost surely and requires some extra assumptions, the  two in M02 and RH99 are not sharp
and none of them is in the ABP format.
There has never been a sharp upper bound for the ABP (or RBP) of the regression median before Proposition 3.2.

\vs
\vs
\noin
{(II)} \tb{Finite sample versus asymptotic breakdown point, the merit of FSBP}.\vs
With the asymptotic breakdown value (or the limit of the finite-sample breakdown) of $\mb{T}^*_{RD}$, $1/3$, had already been given in 
 VAR00 (Theorem 2), RH99 (Theorem 8) and Van Aelst et al. (2002) (Theorem 1), respectively, what is the merit of talking about FSBP?
\vs
\bec
\begin{table}[h!]
\vspace*{-3mm}
\centering
\caption{Average differences between the upper bound of RBP in Proposition 3.2 and the asymptotic breakdown value $1/3$, based on 1000  standard normal samples.}
\vs
\begin{tabular}{c c c c c  c c} 
n& 10& 20& 30& 50& 100& 200\\[-0.0ex]
\hline\\[.1ex]
p=2& 2.447 &5.997  &7.617  &9.105  &10.664 & 11.975 \\[.5ex]
p=3& -7.563 & -2.118&0.387 & 3.072& 5.897 & 8.149 \\[.5ex]
p=5&-20.523 & -13.053 &-9.190 &-5.083 & -0.941 & 1.927 \\[.5ex]
\hline
\end{tabular}
\label{t2}
\end{table}
\vspace*{-8mm}
\enc
This limiting result, 1/3, can  be obtained directly from  (A) of Proposition 3.2 if samples come from  the assumed 
distribution in above references, in this case $k^*(\mb{Z}^n)=n\mbox{RD}(\bs{\beta}^*_{RD}, \mb{Z}^n)$ approaches $ n/2$ as $n\to \infty$ (see Theorems 6 and 7 of RH99). Namely, Proposition 3.2 recovers Theorem 8 of RH99 and Theorem 1 of Van Aelst et al. (2002).
The estimator $\mb{T}^*_{RD}$, however, has to be used in the  finite-sample practice, and the limit $1/3$ is not that informative. For example, the latter implies that
to break down $\mb{T}^*_{RD}$, one must use $n/3$ contaminating points. However, (B) of Proposition 3.2 asserts that one just needs $k^*(\mb{Z}^n)-p+1$ contaminating points.
The differences of  ($k^*(\mb{Z}^n)-p+1)/n$ with $1/3$ in finite-sample cases (especially small sample sizes) are given in Table \ref{t2} or Figure 2 below.
\vs
The table entries reveal once again the merit of the upper bound of the RBP for the FSBP of $\mb{T}^*_{RD}$ because it is
quite different from 
the asymptotic breakdown value (ABV) 1/3 in all cases considered. For example, in p=2 case, the upper bound of RBP indicates that $\mb{T}^*_{RD}$ has a resistance rate for contamination much higher than the ABV 1/3; in $p=3$ case, ABV ranges from overestimating the contamination rate for small $n$s to underestimating the rate for $n>20$;  in $p=5$ case, ABV overestimates the contamination rate for most of $n$s considered.
The ABV, 1/3, is irrelevant for these finite-sample cases because it over-estimates systematically the FSBP of the $\mb{T}^*_{RD}$ in small sample $n$s and large $p$s cases whereas it underestimates the  FSBP for large $n$s with respect to the simulated data, the upper bound of the RBP in Proposition 3.2 increases when $n$ increases for a fixed $p$ and decreases as $p$ increases for a fixed $n$.
 \vs
\bec
\begin{figure}[ht]
\vspace*{-0mm}
    \centering
    \begin{subfigure}[ht]{0.3\textwidth}
        \includegraphics[width=4cm, height=4cm]{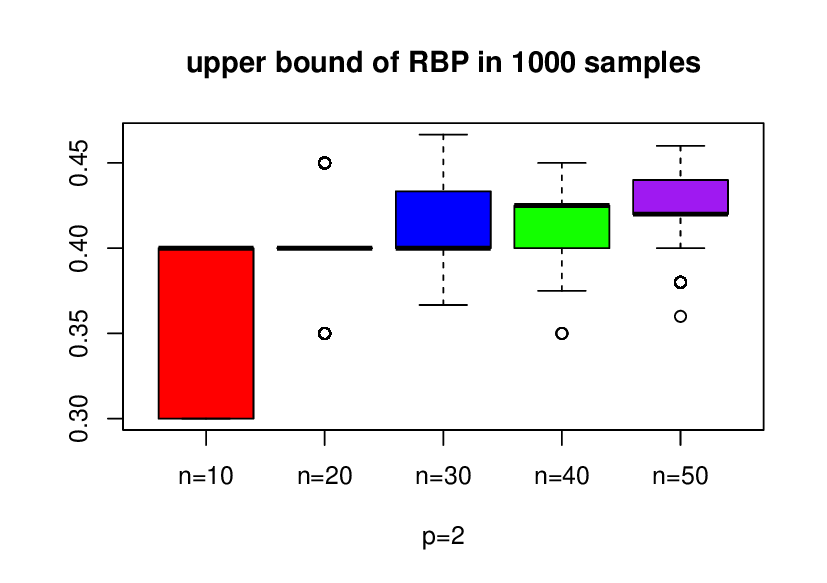}
        \caption{p=2}
        \label{fig:no-camtami}
    \end{subfigure}
     \begin{subfigure}[ht]{0.3\textwidth}
        \includegraphics[width=4cm, height=4cm]{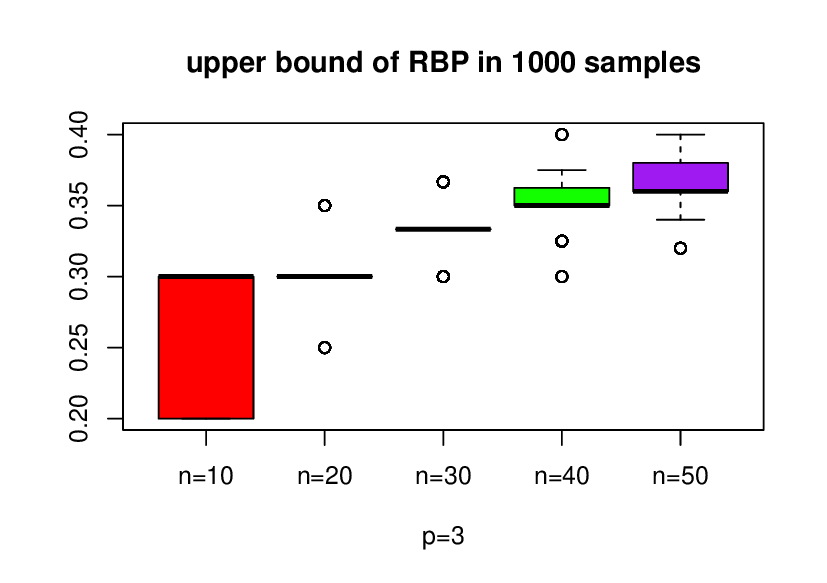}
        \caption{p=3}
        \label{fig:one-contami}
    \end{subfigure}
    \begin{subfigure}[ht]{0.3\textwidth}
        \includegraphics[width=4cm, height=4cm]{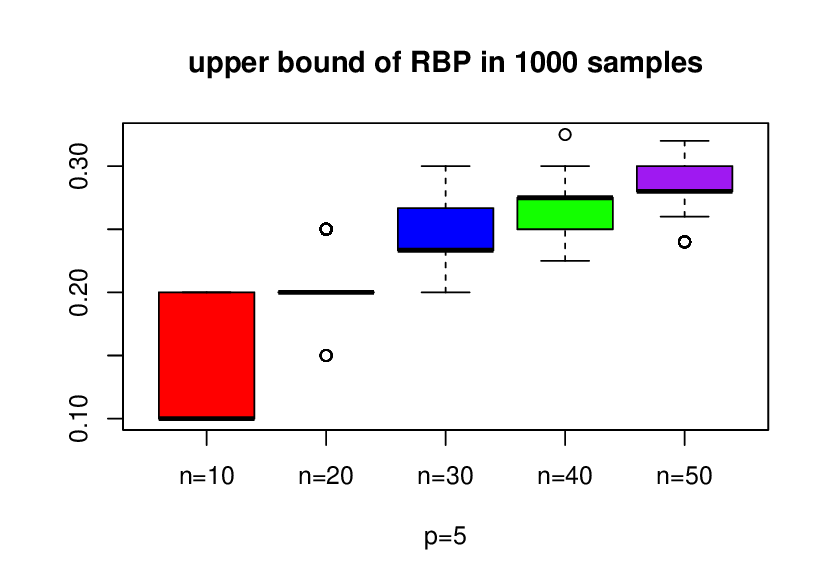}
        \caption{p=5}
        \label{fig:contami}
    \end{subfigure}
    \caption{\small Boxplots for the upper bound of RBP in Proposition 3.2 for $\mb{T}^*_{RD}$  based on 1000  standard normal samples for three  $p$s and five  $n$s.}
    \label{fig:boxplot}
\end{figure}
\vspace*{-8mm}
\enc

Simulation results of the upper bound of RBP in Proposition 3.2 for  $\mb{T}^*_{RD}$ in $1000$ samples can also be displayed graphically in terms of their distributions such as in Figure \ref{fig:boxplot}.
\vs
Inspection of the figure reveals that (i) the upper bound of the RBP in Proposition 3.2 is  always lower than the ABV $1/3$ when $n\leq 20$ and $p>2$,
(ii) it decreases as $p$ increases for a fixed $n$ and increases as $n$ does for a fixed $p$
and (iii) outliers exist in various cases, including $p=2, n= 20, 40, 50$; $p=3, n= 20, 30,40, 50$ and $p=5, n=20,  40, 50$. All these observations and results demonstrate the merit of the FSBP and the relevance of the bounds in Proposition 3.2 (and the irrelevance of the ABV $1/3$) in the finite-sample practice.\vs
\vs

\begin{figure}[t!]
\vspace*{-15mm}
\includegraphics[width=\textwidth]{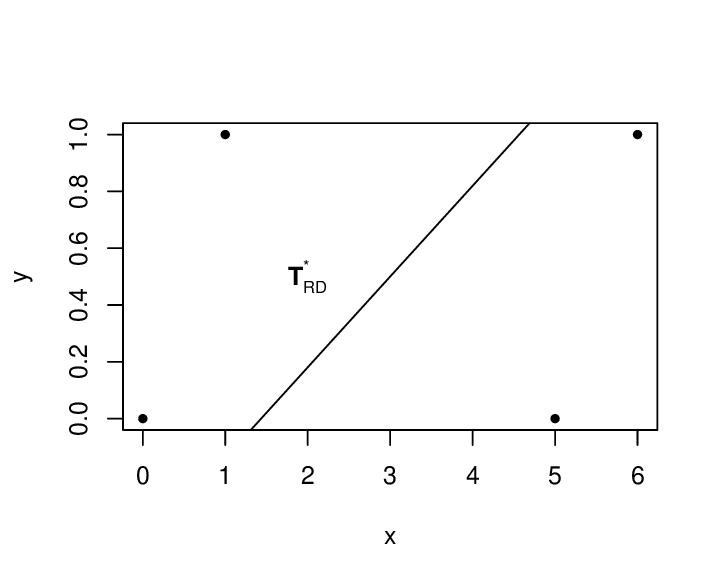}
\vspace*{-13mm}
\caption{\small Four in general position points (represented by four filled circles) located in $\R^2$. Six lines formed, $\{(0, 1), (-5, 1), (1, 0), (0, 0), (0, 1/6), (5/4, -1/4)\}$,  in (intercept, slope) form, each connecting two data points (unillustrated in the figure). 
Each line attains the maximum regression depth $1/2$. However, the average of these deepest lines, $T^*_{RD}$: $(-11/4, ~23/12)/6$ has regression depth 0.}
\label{fig-R2-6-P}
\vspace*{-0mm}
\end{figure}

 \noin
{(III)} \tb{Justification of regression by the maximum depth estimator (median)}.\vs
\noin Proposition 3.2 reveals the intrinsic connection between the breakdown point and the maximum depth value.
This kind of connection was 
 also discussed in M02. 
This intrinsic connection
clearly justifies employment of the maximum depth median as a robust alternative to the traditional regression 
 estimators  since the former is much more robust both in the finite-sample sense and in the asymptotic sense as well.
\vs
\vs
\noin
{(IV)} \tb{Location counterpart and other related results}. \vs
\noin The location counterpart of RD and $\bs{\beta}^*_{RD}$ are respectively halfspace depth (Tukey (1975)) and halfspace median (HM).
The finite-sample breakdown point of the latter has been investigated thoroughly in the literature, e.g., Donoho (1982), Donoho and Gasko (1992) (DG92), Chen (1995) (C95), Chen and Tyler (2000) (CT00),  and Liu {et al}. (2017) (LZW17). \vs In summary, the asymptotic breakdown point of the HM  can be as high as $1/3$ under symmetry and other assumptions (see, C95, CY00, and DG92), also see LZW17 (Proposition 2.10) where only $\mb{X}^n$ is assumed to be IGP. The exact expression of the FSBP of HM is given in LZW17 under two assumptions (i) the $\mb{X}^n$ is IGP and (ii)
a special contaminating scheme: all contaminating points lie at the same site.
It seems that the idea of the proof of Proposition 3.2 could be extended to establish the bounds of the FSBP for HM with an arbitrary (more general) contamination scheme. 
\vs
The exact FSBP of the projection regression depth median, $\mb{T}^*_{PRD}$, a major competitor of the $\mb{T}^*_{RD}$,  has been investigated and established in Zuo (2019). The asymptotic breakdown point of the $\mb{T}^*_{PRD}$ reaches the highest possible value of $50\%$. \vs
\vs
\noin
{(V)} \tb{Computation of regression median}.\vs
\noin The computation of RD and $\mb{T}^*_{RD}$ is challenging and has been  discussed in
RH99 briefly, in Rousseeuw and Struyf (1998), in Van Aelst, Rousseeuw, Hubert, and Struyf (2002),  and in Liu and Zuo (2014). An  \textsf{R} package ``mrfDepth"
has been developed by Segaert, Hubert, Rousseeuw,  Raymaekers, and Vakili (2020). Like most other high breakdown point methods, $\mb{T}^*_{RD}$
has to be computed approximately, which might affect its actual finite-sample breakdown value. 
\vs
\vs
\begin{center}
{\textbf{\large Acknowledgment}}
\end{center}
The author thanks Hanshi Zuo, Hanwen Zuo and Prof. Wei Shao 
for their careful proofreading
and the AE and two anonymous referees for their insightful and constructive comments and suggestions, all of
which have led to improvements in the manuscript.
\vs
\vs
\noindent
\tb{\Large Appendix}
\vs\vs
\noin
\tb{Proof of Proposition 3.2}
\vs
\vs
\noin
{\large{\tb{Part (A)}}}
\vs
\noin
\tb{(i)
 We claim that}
\tb{$\mb{m< k^*(\mb{Z}^n)-p+1}$ contaminating points are not enough to breakdown $\mb{T}^*_{RD}$ in the addition manner} (that is, the RHS of (\ref{lower.eqn}) is a lower bound of ABP for $\mb{T}^*_{RD}$).\vs
Assume, otherwise, that $m< k^*(\mb{Z}^n)-p+1$  contaminating points are enough to breakdown $\mb{T}^*_{RD}$. That is, 
$$\sup_{\mb{Y}^m}\|\mb{T}^*_{RD}(\mb{Z}^n+ \mb{Y}^m)\|=\infty.$$
Denote $\bs{\beta}^*_{RD}(\mb{Z}^n+ \mb{Y}^m)$ (slight notation abuse) as the maximizer of RD at $\mb{Z}^n+ \mb{Y}^m$
which has the maximum norm among all
the RD maximizers. There are at most finitely many RD maximizers for a fixed finite-sample size $n+m$ (assume, without loss of generality (w.l.o.g.), that the hyperplane determined by  a
RD maximizer contains at least $p$ sample points).
\vs
Notice that $\mb{T}^*_{RD}(\mb{Z}^n+ \mb{Y}^m)$ is not necessarily identical to $\bs{\beta}^*_{RD}(\mb{Z}^n+ \mb{Y}^m)$ in the non-unique maximizer case, and the RD of the former could be smaller than the latter and less than the maximum depth value due to the average.  For example,
in the case of the data set $\{(0,0)^{\top}, (1,1)^{\top}, (5,0)^{\top}, (6,1)^{\top}\}$ in $\R^2$, the former is $0$ whereas the latter is $1/2$. See Figure \ref{fig-R2-6-P}.
This is one of the reasons we will treat
$\bs{\beta}^*_{RD}(\mb{Z}^n+ \mb{Y}^m)$ in the sequel instead of  $\mb{T}^*_{RD}(\mb{Z}^n+ \mb{Y}^m)$, in order to compensate for the weakness of the latter.

\bec
\enc

\vspace*{-10mm}
\noin
The other reason, which is the more important one, is
$$
 \sup_{\mb{Y}^m}\|\mb{T}^*_{RD}(\mb{Z}^n+ \mb{Y}^m)\|=\infty \mbox{~if and only if~} \sup_{\mb{Y}^m}\|\bs{\beta}^*_{RD}(\mb{Z}^n+ \mb{Y}^m)\|=\infty.  
$$
The above implies that either
\bi
\item[]\tb{(I)} $|\bs{\beta}^*_{RD}(\mb{Z}^n+ (\mb{Y}^m)_j)_1|\to \infty$ and $\|\bs{\beta}^*_{RD}(\mb{Z}^n+ (\mb{Y}^m)_j)_2\|$ is finite, or\vs
\item[]\tb{(II)}
$\|\bs{\beta}^*_{RD}(\mb{Z}^n+ (\mb{Y}^m)_j)_2\|=\big|\tan(\theta_{\bs{\beta}^*_{RD}(\mb{Z}^n+ (\mb{Y}^m)_j)})\big|\to \infty,$
\ei
along a sequence of $(\mb{Y}^m)_j$ as $j\to \infty$, where the subscripts $1$ and $2$ correspond to the intercept and non-intercept terms, respectively, as in the case $\bs{\beta}=(\beta_1, \bs{\beta}^{\top}_2)^{\top}$.
\vs
\noin
\tb{Case (I)}.\vs
\noin
 Assume that the hyperplane $H_{\bs{\beta}^*_{RD}(\mb{Z}^n+ (\mb{Y}^m)_j)}$ intersects the horizontal hyperplane $y=0$ at the hyperline
$l_v(\bs{\beta}^*_{RD}(\mb{Z}^n+ (\mb{Y}^m)_j)$ (when the two do not intersect then we assume that the hyperline exists at infinity, the arguments hereafter go through). Since the intercept term of ${\bs{\beta}^*_{RD}(\mb{Z}^n+ (\mb{Y}^m)_j)}$ approaches  infinity and the $\|\bs{\beta}^*_{RD}(\mb{Z}^n+ (\mb{Y}^m)_j)_2\|$ is finite, the hyperplane  no longer contains  any points from $\mb{Z}^n$.
There are at most $m$ contaminating points from $\mb({Y}^m)_j$ on the hyperplane. Therefore it is readily seen that
\begin{align}
m+(p-1)\geq m &\geq (m+n)\mbox{RD}(\bs{\beta}^*_{RD}(\mb{Z}^n+ (\mb{Y}^m)_j), \mb{Z}^n+ (\mb{Y}^m)_j)\nonumber \\[1ex]
&=k^*(\mb{Z}^n+ (\mb{Y}^m)_j)\geq k^*(\mb{Z}^n), \label{inequality-0.eqn}
\end{align}
where the first inequality is trivial and the second inequality follows from the facts that (i) one can tilt $H_{\bs{\beta}^*_{RD}(\mb{Z}^n+ (\mb{Y}^m)_j)}$ along $l_v(\bs{\beta}^*_{RD}(\mb{Z}^n+ (\mb{Y}^m)_j)$ to a vertical position without touching any points from $\mb{Z}^n$ and (ii) the definition of RD or  lemma 3.1. The third equality follows from (\ref{k*.eqn}) and the definition of $\bs{\beta}^*_{RD}(\mb{Z}^n+ (\mb{Y}^m)_j)$ above. The last inequality follows from the fact that $(m+n)\mbox{RD}(\bs{\beta}, \mb{Z}^n+ \mb{Y}^m)\geq n\mbox{RD}(\bs{\beta}, \mb{Z}^n)$ in light of (\ref{rd-proof.eqn}).
\vs

\noin
\tb{Case (II)}. \vs
\noin
\tb{If, there exists a finite $\mb{j}$ such that $\mb{{\theta}_{\bs{\beta}^*_{RD}(\mb{Z}^n+ (\mb{Y}^m)_j} =\bs{\pi}/2}$,}
then at most  m contaminating points from $(Y^m)_j$ are on the hyperplane $H_{\bs{\beta}^*_{RD}(\mb{Z}^n+ (\mb{Y}^m)_j)}$ which contains  at most $p-1$ points from $Z^n$. The latter is due to the fact that $Z^n$ is IGP and the intersection hyperline between  $H_{\bs{\beta}^*_{RD}(\mb{Z}^n+ (\mb{Y}^m)_j)}$ and the horizontal hyperplane $y=0$ is a $(p-2)$ dimensional subspace
of the $p$ dimensional space $(\mb{x}^{\top}, y)$. It is not hard to see that
\begin{align}
m+(p-1)&\geq (m+n)\mbox{RD}(\bs{\beta}^*_{RD}(\mb{Z}^n+ (\mb{Y}^m)_j), \mb{Z}^n+ (\mb{Y}^m)_j)\nonumber \\[1ex]
&=k^*(\mb{Z}^n+ (\mb{Y}^m)_j)\geq k^*(\mb{Z}^n), \label{inequality-1.eqn}
\end{align}
where the first inequality follows from the fact that the vertical hyperplane contains at most $m+(p-1)$ points from $\mb{Z}^n+ (\mb{Y}^m)_j$.
The second equality follows from (\ref{k*.eqn}) and the definition of $\bs{\beta}^*_{RD}(\mb{Z}^n+ (\mb{Y}^m)_j)$ above. The last inequality follows from the fact that $(m+n)\mbox{RD}(\bs{\beta}, \mb{Z}^n+ \mb{Y}^m)\geq n\mbox{RD}(\bs{\beta}, \mb{Z}^n)$ in light of (\ref{rd-proof.eqn}).
\vs
\noin
\tb{Otherwise (i.e. $\mb{\theta_{\bs{\beta}^*_{RD}(\mb{Z}^n+ (\mb{Y}^m)_j)} <\bs{\pi}/2}$, for any $\mb{j}$ and $\mb{\theta_{\bs{\beta}^*_{RD}(\mb{Z}^n+ (\mb{Y}^m)_j)} \rightarrow \bs{\pi}/2}$),}
assume that $x_{11}, \cdots, x_{n1}$ are the  $n$ first coordinates of $\mb{Z}_1,\cdots, \mb{Z}_n$ (note that $\mb{Z}_i=(\mb{x}^{\top}_i, y_i)^{\top}=(x_{i1},\cdots, x_{i(p-1)}, y_i)^{\top}$). Now order $x_{11}, \cdots, x_{n1}$
and assume that $x_{i1}$ and $x_{j1}$ are two distinct and consecutive (after ordering) first coordinates and let $\delta=x_{j1}-x_{i1}>0$.   
Consider, a vertical hyperplane
$H_v=\{ (\mb{x}^{\top}, y)^{\top}\in\R^p, x_1=(x_{i1}+x_{j1})/2 \}$ that contains no data points from $\mb{Z}^n$. Assume, w.l.o.g., for any large enough $j$ and $(\mb{Y}^{m})_j$ and $H_{\bs{\beta}^*_{RD}(\mb{Z}^n+ (\mb{Y}^m)_j)}$, $H_v$ intersects with  $H_{\bs{\beta}^*_{RD}(\mb{Z}^n+ (\mb{Y}^m)_j)}$ at $l_v(\bs{\beta}^*_{RD}(\mb{Z}^n+ (\mb{Y}^m)_j)$.
\vs
Clearly,
there exists a narrow vertical hyperstrip centered at $H_v$ (with its two boundary hyperplanes parallel to $H_v$, e.g., $H_{v1}=\{ (\mb{x}^{\top}, y)^{\top}\in\R^p, x_1=(x_{i1}+x_{j1})/4 \}$ and $H_{v2}=\{ (\mb{x}^{\top}, y)^{\top}\in\R^p, x_1=3(x_{i1}+x_{j1})/4 \}$ as its two boundary hyperplanes), and within the hyperstrip/hyperslab there are no data points from $\mb{Z}^n$.
\vs
Now, when one tilts the hyperplane $H_{\bs{\beta}^*_{RD}(\mb{Z}^n+ (\mb{Y}^m)_j)}$ (which is already almost vertical for a large enough $j$) along $l_v (\bs{\beta}^*_{RD}(\mb{Z}^n+ (\mb{Y}^m)_j)$ to its eventual vertical position of $H_v$, it is readily apparent that
\begin{align}
m+p-1&\geq (m+n) \min_{fr}(l_v(\bs{\beta}^*_{RD}(\mb{Z}^n+ (\mb{Y}^m)_j), \mbox{Pr}_{n+m})\nonumber\\[1ex]
&\geq \mbox{RD}(\bs{\beta}^*_{RD}(\mb{Z}^n+ (\mb{Y}^m)_j)= k^*(\mb{Z}^n+ (\mb{Y}^m)_j)\nonumber\\[1ex]
&\geq k^*(\mb{Z}^n), \label{inequality-2.eqn}
\end{align}
 where $\mbox{Pr}_{m+n}$ stands for the empirical distribution based on $\mb{Z}^n+ (\mb{Y}^m)_j$.
 The first inequality follows from (i) the definition of $\min_{fr}(l_v(\bs{\beta}, \mbox{Pr}_n)$ and (ii) the fact that there is one way of tilting the hyperplane $H_{\bs{\beta}^*_{RD}(\mb{Z}^n+ (\mb{Y}^m)_j)}$ along $l_v (\bs{\beta}^*_{RD}(\mb{Z}^n+ (\mb{Y}^m)_j))$ so that no original data points from $\mb{Z}^n$ (except at most $p-1$ points from $\mb{Z}^n$ that are already on the hyperpline $l_v$) are touched during the movement and the points it can touch are at most all the m contaminating points (since for larger enough $j$, $H_{\bs{\beta}^*_{RD}(\mb{Z}^n+ (\mb{Y}^m)_j)}$ is within the  hyperstrip/hyperslab formed by $H_{v1}$ and $H_{v2}$).
 \vs
  The second inequality above follows  (\ref{rd-proof.eqn}) and the third equality follows from the introduction of $\bs{\beta}^*_{RD}(\mb{Z}^n+ \mb{Y}^m)$ at the beginning and  from  (\ref{k*.eqn}). Finally, the fourth inequality comes from the fact that $(m+n)\mbox{RD}(\bs{\beta}, \mb{Z}^n+ \mb{Y}^m)\geq n\mbox{RD}(\bs{\beta}, \mb{Z}^n)$ in light of (\ref{rd-proof.eqn}).
\vs
All the inequalities in (\ref{inequality-0.eqn}), (\ref{inequality-1.eqn}) and (\ref{inequality-2.eqn})  lead to a contradiction.
Thus, $m< k^*(\mb{Z}^n)-p+1$ contaminating points are not enough to break down $\mb{T}^*_{RD}$.
\vs
\noindent
\tb{(ii) We claim that} \tb{m=$k^*(\mb{Z}^n)-p+1$ points are enough to break down $\mb{T}^*_{RD}$ in the addition manner} (that is, the RHS of (\ref{lower.eqn}) is also an upper bound of ABP for $\mb{T}^*_{RD}$).
\vs
Let $l_h$ be a hyperline (with dimension $p-2$) in the (p-1)-dimensional $\mb{x}$ space ($y=0$, the subspace of $(\mb{x}^{\top}, y)^{\top}$-space in $\R^p$) that contains  $p-1$ points of $\mb{x}^{\top}_is$ 
 and $H_{v*}$ be the corresponding vertical hyperplane that intercepts with the horizontal hyperplane $y=0$ at the hyperline $l_h$.
\vs
Construct another hyperplane $H_{\bs{\beta}_o}$ that also intercepts with the horizontal hyperplane $y=0$ at the hyperline $l_h$ and is almost vertical.
Place $m$ contaminating points in $\mb{Y}^m$ at a point $\mb{Z}=(\mb{x}^{\top}, y)^{\top}$ on $H_{\bs{\beta}_o}$, where $y$ could be arbitrarily large (defined below) and the position of $\mb{Z}$ is described below.
Denote 
the resulting hyperplane based on $H_{\bs{\beta}_o}$ by $H_{\bs{\beta}_c}$
and the resulting 
data set by $\mb{Z}^{n+m}:=\mb{Z}^n+ \mb{Y}^m$,
It is readily seen that $(n+m) \mbox{RD}(\bs{\beta}_c, \mb{Z}^{n+m})\geq k^*(\mb{Z}^n)$ in light of (\ref{rd-proof.eqn}).

\vs
If we can show that $\bs{\beta}_c$ attains the maximum RD with respect to (w.r.t.) $\mb{Z}^{n+m}$, then when we tilt $H_{\bs{\beta}_c}$ to the vertical position $H_{v*}$,  we have added $m$ points to break down $\mb{T}^*_{RD}$, which is the average of all maximizers of RD w.r.t. $\mb{Z}^{n+m}$. 
Equivalently, we have to show that for any given $\bs{\beta} \in \R^p$, 
RD$(\bs{\beta}, \mb{Z}^{n+m})\leq k^*(\mb{Z}^n)/(n+m)$.
\vs
Denote by $H_{\bs{\beta}}$ the unique hyperplane determined by $\bs{\beta}$ (through $y=(1,\mb{x}^{\top})\bs{\beta}$). Consider two cases in the sequel: (a) $\mb{Z}\in H_{\bs{\beta}}$ (b) $\mb{Z}\not\in H_{\bs{\beta}}$.

\bec
\begin{figure}[t!]
\vspace*{-20mm}
\includegraphics[width=\textwidth]{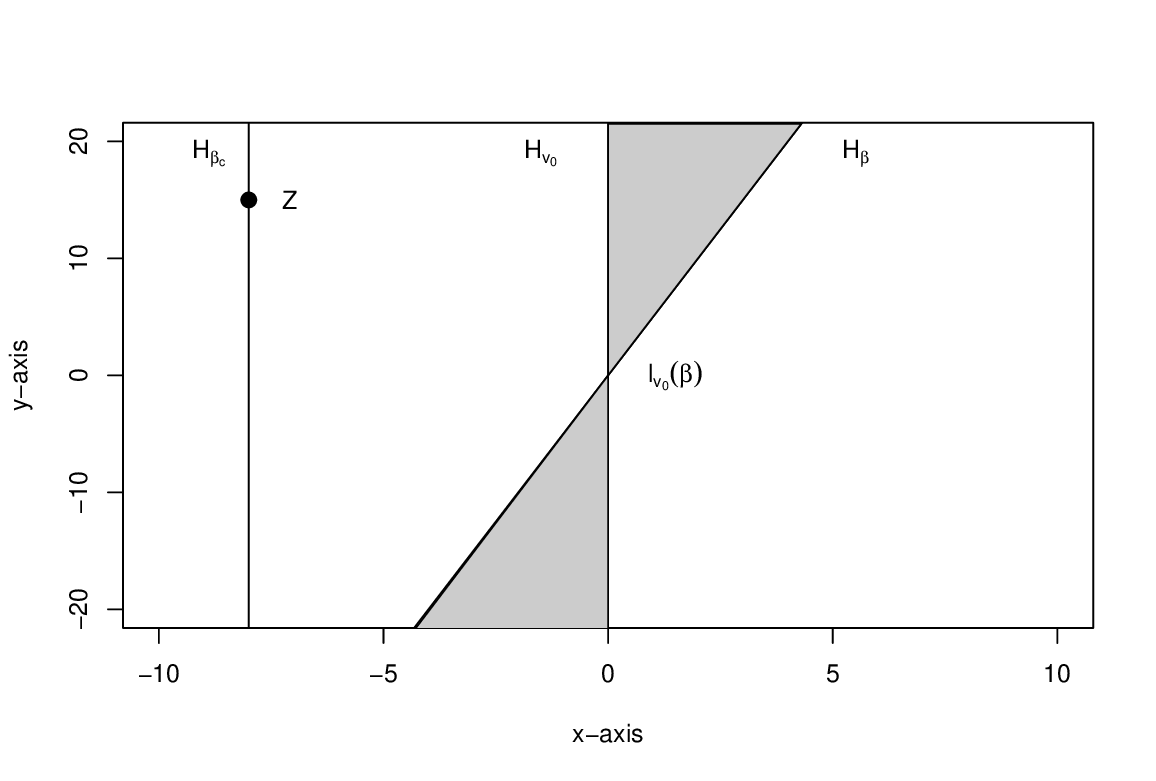}
\vspace*{-13mm}
 \caption{\small A two-dimensional vertical cross-section of a figure in  $\R^p$. Assume that tilting $H_{\bs{\beta}}$ to a vertical position $H_{v_0}$ along hyperline $l_{v_0}(\bs{\beta})$ crossing the shaded double wedge touches the minimum fraction of data points in $\mb{Z}^n$ among all possible $l_v$s. That is, the fraction is exactly RD$(\bs{\beta}, \mb{Z}^n)$. Note that the $y$ component of $Z$ controls its vertical position and $\mb{x}$ component of Z can control the horizontal position of $H_{\bs{\beta_c}}$.
  }
 \label{fig-4-proof-3}
\end{figure}
\enc
\vspace*{-8mm}
\noin
\tb{(a) $\mb{Z}\in H_{\bs{\beta}}$}. Assume that $H_{\bs{\beta}}$ intersects with the horizontal hyperplane ($y=0$) at hyperline $l(\bs{\beta})$. Call the vertical hyperplane that intercepts $y=0$ at $l(\bs{\beta})$ as $H_v$.
Apparently, $l(\bs{\beta})$ contains at most $p-1$ points from $\mb{Z}^n$. Construct a narrow vertical hyperstrip centered at $l(\bs{\beta})$ with its two parallel boundary hyperplanes parallel to $H_{v}$ (define $\delta=\min_{\mb{x}_i\not \in l(\bs{\beta}),~ i\in\{1,\cdots, n\}}d(\mb{x}_i, l(\bs{\beta}))$, where $d (\cdot, \cdot)$ stands for the distance of a point to a set. Introduce two hyperplanes that are parallel to $H_{v}$ and with distance $\delta/2$ to $H_v$), within the hyperstrip there are no data points from $\mb{Z}^n$ except those with $\mb{x}$-component on $l(\bs{\beta})$. By choosing large enough $y$ so that
 $H_{\bs{\beta}}$ is almost vertical and is within the  vertical hyperstrip. 
Now tilting $H_{\bs{\beta}}$ to the position of $H_v$, it touches at most $m+(p-1)$ points from $\mb{Z}^{n+m}$ (at most $p-1$ points with $\mb{x}$-component  on the $l(\bs{\beta})$ and  $m$ contaminating points at $\mb{Z}$), that is, RD$(\bs{\beta}, \mb{Z}^{n+m})\leq k^*(\mb{Z}^n)/(n+m)$. (this holds true obviously when $H_{\bs{\beta}}$ is parallel to the horizontal hyperplane $y=0$).
\vs
\noin
\tb{(b) $\mb{Z}\not\in H_{\bs{\beta}}$}.  Assume, w.l.o.g., that $H_{\bs{\beta}}$ contains $p$ points from $\mb{Z}^n$ (after all we only care about $H_{\bs{\beta}}$s that contain  most points from $\mb{Z}^n$ so that $\bs{\beta}$ can have as large as possible RD value w.r.t. $\mb{Z}^{n+m}$ that is no greater than $k^*(\mb{Z}^n)/(n+m)$).
Let $H_{v_0}$ be the vertical hyperplane that intersects with the horizontal hyperplane ($y=0$) at $l_{v_0}$ such that when tilting $H_{\bs{\beta}}$ along $l_{v_0}$ to $H_{v_0}$ in one of two ways (assume, w.l.o.g., that it is counter-clockwise, see Figure \ref{fig-4-proof-3}), the number of points in $\mb{Z}^n$ touched is exactly $k:=n\mbox{RD}(\bs{\beta}, \mb{Z}^n)\leq k^*(\mb{Z}^n)$.
\vs
Note that there are only \emph{finitely many}  $H_{\bs{\beta}}$s considered above. Put all $\bs{\beta}$s that with counter-clockwise tilting above and touched
$k$ points of $\mb{Z}^n$ into a group called $G_{c-cw}$ . Those $\bs{\beta}$s that with clockwise tilting above and touched
$k$ points of $\mb{Z}^n$  into a group called $G_{cw}$.
\vs
Now we show that for any $\bs{\beta}$ in $G_{c-cw}$, RD$(\bs{\beta}, \mb{Z}^{n+m})\leq k^*(\mb{Z}^n)/(n+m)$. (Treatments for $\bs{\beta}$ in $G_{cw}$ are similar and thus are skipped).\vs

Tilting $H_{\bs{\beta}}$ along $l_{v_0}$ to vertical position $H_{v_0}$, call the region on $H_{\bs{\beta}_c}$ touched during the tilting process as $R(\bs{\beta}, l_{v_0}, \bs{\beta}_c)$. Since there are finitely many  $H_{\bs{\beta}}$s hence finitely many $R(\bs{\beta}, l_{v_0},\bs{\beta}_c)$s. Let $y_{max}$ be the maximum of finitely many $y$s, where $y$ is the supremum of $y$-components of (points in) each region $R(\bs{\beta}, l_{v_0}, \bs{\beta}_c)$.
Consequently, there is at least one point $\mb{Z}$ (with $y$ component greater than $y_{max}$) on $H_{\bs{\beta}_c}$ that does not lie in the union of all $R(\bs{\beta}, l_{v_0},\bs{\beta}_c)$s.
Place $m$ contaminating points at $\mb{Z}$ so that when tilting $H_{\bs{\beta}}$ along $l_{v_0}$ counter-clockwise to $H_{v_0}$, it does not touch
$\mb{Z}$, see Figure \ref{fig-4-proof-3}.
\vs
Tilting $H_{\bs{\beta}}$, counter-clockwise along $l_{v_0}$ to $H_{v_0}$, the number of the data points in $\mb{Z}^{n+m}$ touched is $k$ (since only points in $\mb{Z}^n$ originally lie in the two wedges passed might be touched).
On the other hand, when tilting it in the clockwise way, the number of points touched in $\mb{Z}^{n+m}$,  denoted by $q$,  is  at least $p+m=k^*(\mb{Z}^n)+1$. Now we have
\begin{align}
\mbox{RD}(\bs{\beta}, \mb{Z}^{n+m})&=\inf_{l_v(\bs{\beta})}\min_{fr}(l_v(\bs{\beta}), \mb{Z}^{n+m})\leq \min_{fr}(l_{v_0}(\bs{\beta}), \mb{Z}^{n+m})\\
&= \min\{k, q\}/(n+m)=k/(n+m)\leq k^*(\mb{Z}^{n})/(n+m),
\end{align}
where the first equality follows from Lemma 3.1, the second inequality is trivial, the third equality follows from above discussions on the number of points touched by tilting $H_{\bs{\beta}}$ in two ways and  from the definition of $\min_{fr}$ given in (or before) Lemma 3.1, the fourth equality is trivial, so is the last inequality. 
We complete the proof of part \tb(A). 

\vs
\noindent
{{\large{\tb{Part (B)}}}\vs
\noindent
\tb{We claim that} \tb{m=$k^*(\mb{Z}^n)-p+1$ points are enough to break down $\mb{T}^*_{RD}$ in the replacement manner.}
\vs
The proof of this part is  similar to that of (ii) above and details are skipped.
\hfill\pend

\vs

\end{document}